\documentclass[12pt,dmathesis]{article}
\usepackage{amstext}
\usepackage[fleqn]{amsmath}
\usepackage{amsfonts}
\usepackage{amssymb}
\usepackage{amsthm}
\usepackage{newlfont}
\usepackage{graphicx}
\usepackage{sectsty}
\theoremstyle{plain}
\theoremstyle{theorem}
\newtheorem{thm}{Theorem}[section]
\theoremstyle{example}

\theoremstyle{corollary}
\newtheorem{cor}{Corollary}[section]
\theoremstyle{lemma}

\theoremstyle{proposition}

\theoremstyle{axiom}

\theoremstyle{notation}

\theoremstyle{fact}

\theoremstyle{definition}
\newtheorem{defn}{Definition}[section]
\theoremstyle{remark}
\newtheorem{rem}{Remark}[section]
\numberwithin{equation}{section}
\begin{document}
\title{Certain new formulas for bibasic Humbert hypergeometric functions $\Psi_{1}$ and $\Psi_{2}$}
\author{Ayman Shehata \thanks{%
E-mail:aymanshehata@science.aun.edu.eg, drshehata2006@yahoo.com, A.Ahmed@qu.edu.sa}\\
{\small Department of Mathematics, Faculty of Science, Assiut University, Assiut 71516, Egypt.}\\
{\small Department of Mathematics, College of Science and Arts, Unaizah 56264, Qassim University,}\\
{\small Qassim, Saudi Arabia.}}
\date{}
\maketitle{}
\begin{abstract}
The main aim of the present work is to give some interesting the $q$-analogues of various $q$-recurrence relations, $q$-recursion formulas, $q$-partial derivative relations, $q$-integral representations, transformation and summation formulas for bibasic Humbert hypergeometric functions $\Psi_{1}$ and $\Psi_{2}$ on two independent bases $q$ and $p$ of two variables and some developments formulae, believed to be new, by using the conception of $q$-calculus. Finally, some interesting special cases and straightforward identities connected with  bibasic Humbert hypergeometric series of the types $\Psi_{1}$ and $\Psi_{2}$ are established when the two independent bases $q$ and $p$ are equal.
\end{abstract}
\section{Introduction}
Quantum calculus is the modern name for the investigation of calculus without limits. The quantum calculus or $q$-calculus began with Jackson in the
early twentieth century, but this kind of calculus had already been worked out by Euler and Jacobi. $q$-calculus appeared as a connection between mathematics and physics. It has a lot of their applications in many fields of mathematics areas such as number theory, an engineering, combinatorics, orthogonal polynomials, basic hypergeometric functions,  quantum theory, physics, mechanics, the theory of relativity and other sciences, see for example \cite{fn, ljw, py, si}. Recently, there have been many studies on $q$-calculus. Andrews \cite{an}, Bytev and Zhang \cite{bz}, Verma \cite{vs}, Verma and Sahai \cite{vs1}, Verma and Sarasvati \cite{vsy} and Yadav et al. \cite{ypv} discussed recursion formulas and transformations for $q$-hypergeometric series. Sahai and Verma \cite{vs2}, Sears \cite{se}, Srivastava \cite{sr1, sr4} and  Upadhyay \cite{up} derived transformations and summation formulas for bilateral basic hypergeometric series. In particular, Jackson \cite{j1, j2} was the first to study basic Appell series. Agarwal has developed some properties of basic Appell series \cite{ajc1}, and Slater \cite{sl} applies contour integral techniques to such series. In a recent paper the authors \cite{gim, ga1, ga2, gr1, sr2, sr3, sr5} have developed very general transformations involving bibasic $q$-Appell hypergeometric series on two unconnected bases ($q$ and $p$, $0<|q|<1$, $0<|p|<1$, $q,p\in\mathbb{C}$). Such generalized basic Appell hypergeometric series were called bibasic. Earlier, the author in \cite{sh1, sh2, sh3, sh4} have developed and studied some relations of the $(p,q)$-Bessel, $(p,q)$-Humbert functions and basic Horn functions $H_{3}$, $H_{4}$, $H_{6}$ and $H_{7}$. Motivated by in the previous work \cite{sh5, sh6}, many terminating summation and transformation formulas for bibasic $q$-Humbert hypergeometric series are derived and proved by using contiguous relations which extend most of the results due to Shehata and have developed their transformation theory. The paper is concluded by obtaining some $q$-recurrence relations, $q$-derivatives formulas, $q$-partial derivative relations, $q$-derivatives with respect to the parameters, $q$-integral representations, transformations and summation formulas for these bibasic Humbert hypergeometric functions $\Psi_{1}$ and $\Psi_{2}$ with different bases $p$ and $q$. We think these results are not found in the literature to discuss further consequences of our extensions as special cases.
\subsection{Notations and preliminaries}
First, we start by remembering some elementary definitions and notations of $q$-analogue with $q$-derivative then definition and properties of $q$-integral used (needed) in the $q$-theory. Throughout this study, unless otherwise is stated, the bases $q$ and $p$ will be assumed to be such that $0<|q|<1$, $0<|p|<1$, and $q,p\in\mathbb{C}$, for definiteness. We use $\mathbb{C}$ to denote the set of complex numbers and $\mathbb{N}$ to denote the set of positive integers.

For $q\in\mathbb{C}$ and $0<|q|<1$, the $q$-shifted factorials $(q^{a};q)_{k}$ are defined as
\begin{eqnarray}
\begin{split}
&(q^{a};q)_{k}=\left\{
            \begin{array}{ll}
              \prod_{r=0}^{k-1}(1-q^{a+r}), & \hbox{$k\geq1$;} \\
              1, & \hbox{$k=0$.}
            \end{array}
          \right.\\
&          =\left\{
  \begin{array}{ll}
    (1-q^{a})(1-q^{a+1})\ldots(1-q^{a+k-1}), & \hbox{$k\in\Bbb{N}, a\in\mathbb{C}\setminus \{0, -1, -2,\ldots,1-k\}$;} \\
    1, & \hbox{$k=0,a\in\mathbb{C}$.}
  \end{array}
\right.\label{1.1}
\end{split}
\end{eqnarray}
\begin{defn}
The $q$-integer number is defined as
\begin{equation}
\begin{split}
[\chi]_{q}=\frac{1-q^{\chi}}{1-q},\chi\in \Bbb{N}_{0}.\label{1.2}
\end{split}
\end{equation}
\end{defn}
\begin{defn}
For $a$, $b$, $c \in \mathbb{C}$, $c\neq 0,-1,-2,...$ and $0<|q|<1$, $q\in\mathbb{C}$, the basic hypergeometric series with base $q$ is defined as (see \cite{am,kc, kk, kls, ks})
\begin{equation}
\begin{split}
\;_{2}\phi_{1}(q^{a},q^{b};q^{c};q,x)=\sum_{k=0}^{\infty}\frac{(q^{a};q)_{k}(q^{b};q)_{k}}{(q^{c};q)_{k}(q;q)_{k}}x^{k},\label{1.3}
\end{split}
\end{equation}
for $|x|<1$ and by analytic continuation for other $x\in \mathbb{C}$.
\end{defn}
\begin{defn}
Let $f$ be a function defined on a subset of the complex or real plane. The $q$-difference operator $\mathbb{D}_{x,q}$ is defined \cite{j3} as follows
\begin{equation}
\begin{split}
\mathbb{D}_{x,q}f(x)=\frac{f(x)-f(qx)}{(1-q)x},x\neq0.\label{1.4}
\end{split}
\end{equation}
\end{defn}
\begin{defn}
For $0<|p|<1$, $0<|q|<1$, $p, q\in\mathbb{C}$, we define the bibasic Humbert functions $\Psi_{1}$ and $\Psi_{2}$ on two independent bases $p$ and $q$ as follows
\begin{eqnarray}
\begin{split}
\Psi_{1}(q^{a},p^{b};p^{c},q^{d};q,p,x,y)=&\sum_{\ell,k=0}^{\infty}\frac{(q^{a};q)_{k+\ell}(p^{b};p)_{\ell}}{(p^{c};p)_{\ell}(q^{d};q)_{k}(q;q)_{k}(p;p)_{\ell}}x^{k}y^{\ell},\\
&(p^{c},q^{d}\neq 1, q^{-1}, q^{-2},\ldots,|x|,|y|<1)\label{1.5}
\end{split}
\end{eqnarray}
and
\begin{eqnarray}
\begin{split}
\Psi_{2}(q^{a};p^{b},q^{c};q,p,x,y)=&\sum_{\ell,k=0}^{\infty}\frac{(q^{a};q)_{k+\ell}}{(p^{b};p)_{\ell}(q^{c};q)_{k}(q;q)_{k}(p;p)_{\ell}}x^{k}y^{\ell},\\
&(q^{c}, p^{b}\neq 1, q^{-1}, q^{-2},\ldots,|x|,|y|<1).\label{1.6}
\end{split}
\end{eqnarray}
\end{defn}
\section{Main Results}
In this section, we derive the $q$-analogues and extensions of bibasic Humbert functions $\Psi_{1}$ and $\Psi_{2}$ on two independent bases $q$ and $p$ with their several interesting properties on the same order.
\begin{thm} The functions $\Psi_{1}$ and $\Psi_{2}$ satisfy the following recurrence relations.
\begin{eqnarray}
\begin{split}
&\Psi_{1}(q^{a+1},p^{b};p^{c},q^{d};q,p,x,y)=\Psi_{1}(q^{a},p^{b};p^{c},q^{d};q,p,x,y)+\frac{q^{a}x}{1-q^{d}}\\
&\times \Psi_{1}(q^{a+1},p^{b};p^{c},q^{d+1};q,p,x,y)+\frac{q^{a}}{1-q^{a}}\Psi_{1}(q^{a},p^{b};p^{c},q^{d};q,p,q x,y)\\
&-\frac{q^{a}}{1-q^{a}}\Psi_{1}(q^{a},p^{b};p^{c},q^{d};q,p,q x,q y),q^{a},q^{d}\neq 1,\label{2.1}
\end{split}
\end{eqnarray}
\begin{eqnarray}
\begin{split}
&\Psi_{1}(q^{a+1},p^{b};p^{c},q^{d};q,p,x,y)=\frac{1}{1-q^{a}}\Psi_{1}(q^{a},p^{b};p^{c},q^{d};q,p,x,y)-\frac{q^{a}}{1-q^{a}}\\
&\times \Psi_{1}(q^{a},p^{b};p^{c},q^{d};q,p,x,q y)+\frac{q^{a}x}{1-q^{d}}\Psi_{1}(q^{a+1},p^{b};p^{c},q^{d+1};q,p,x,q y),q^{a},q^{d}\neq 1,\\
&\Psi_{1}(q^{a},p^{b};p^{c},q^{d-1};q,p,x,y)=\Psi_{1}(q^{a},p^{b};p^{c},q^{d};q,p,x,y)\\
&+\frac{q^{d-1}(1-q^{a})x}{(1-q^{d-1})(1-q^{d})}\Psi_{1}(q^{a+1},p^{b};p^{c},q^{d+1};q,p,x,y),q^{d},q^{d-1}\neq 1\label{2.2}
\end{split}
\end{eqnarray}
and
\begin{eqnarray}
\begin{split}
&\Psi_{2}(q^{a+1};p^{b},q^{c};q,p,x,y)=\Psi_{2}(q^{a};p^{b},q^{c};q,p,x,y)+\frac{q^{a}x}{1-q^{c}}\Psi_{2}(q^{a+1};p^{b},q^{c+1};q,p,x,y)\\
&+\frac{q^{a}}{1-q^{a}}\Psi_{2}(q^{a};p^{b},q^{c};q,p,q x,y)-\frac{q^{a}}{1-q^{a}}\Psi_{2}(q^{a};p^{b},q^{c};q,p,q x,q y),q^{a},q^{c}\neq 1,\\
&\Psi_{2}(q^{a+1};p^{b},q^{c};q,p,x,y)=\frac{1}{1-q^{a}}\Psi_{2}(q^{a};p^{b},q^{c};q,p,x,y)-\frac{q^{a}}{1-q^{a}}\\
&\times\Psi_{2}(q^{a};p^{b},q^{c};q,p,x,q y)+\frac{q^{a}x}{1-q^{c}}\Psi_{2}(q^{a+1};p^{b},q^{c+1};q,p,x,q y),q^{a},q^{c}\neq 1,\\
&\Psi_{2}(q^{a};p^{b},q^{c-1};q,p,x,y)=\Psi_{2}(q^{a};p^{b},q^{c};q,p,x,y)\\
&+\frac{q^{c-1}(1-q^{a})x}{(1-q^{c-1})(1-q^{c})}\Psi_{2}(q^{a+1};p^{b},q^{c+1};q,p,x,y),q^{c},q^{c-1}\neq 1.\label{2.3}
\end{split}
\end{eqnarray}
\end{thm}
\begin{proof} To prove the identity (\ref{2.1}). Using the relations
\begin{eqnarray*}
\begin{split}
(q^{a+1};q)_{k+\ell}=&\frac{1-q^{a+k+\ell}}{1-q^{a}}(q^{a};q)_{k+\ell}=\bigg{[}1+q^{a}\frac{1-q^{k+\ell}}{1-q^{a}}\bigg{]}(q^{a};q)_{k+\ell},\\
(q^{a};q)_{k+\ell+1}=&(1-q^{a})(q^{a+1};q)_{k+\ell},
\end{split}
\end{eqnarray*}
and
\begin{eqnarray*}
\begin{split}
1-q^{k+\ell}=&1-q^{k}+q^{k}-q^{k+\ell},\\
=&1-q^{\ell}+q^{\ell}(1-q^{k}),
\end{split}
\end{eqnarray*}
we have
\begin{eqnarray*}
\begin{split}
&\Psi_{1}(q^{a+1},p^{b};p^{c},q^{d};q,p,x,y)-\Psi_{1}(q^{a},p^{b};p^{c},q^{d};q,p,x,y)\\
&=\frac{q^{a}}{1-q^{a}}\sum_{\ell=0,k=1}^{\infty}\frac{(q^{a};q)_{k+\ell}(p^{b};p)_{\ell}}{(p^{c};p)_{\ell}(q^{d};q)_{k}(q;q)_{k-1}(p;p)_{\ell}}x^{k}y^{\ell}+\frac{q^{a}}{1-q^{a}}\sum_{\ell,k=0}^{\infty}\frac{q^{k}(q^{a};q)_{k+\ell}(p^{b};p)_{\ell}}{(p^{c};p)_{\ell}(q^{d};q)_{k}(q;q)_{k}(p;p)_{\ell}}x^{k}y^{\ell}\\
&-\frac{q^{a}}{1-q^{a}}\sum_{\ell,k=0}^{\infty}\frac{q^{k+\ell}(q^{a};q)_{k+\ell}(p^{b};p)_{\ell}}{(p^{c};p)_{\ell}(q^{d};q)_{k}(q;q)_{k}(p;p)_{\ell}}x^{k}y^{\ell}\\
&=\frac{q^{a}}{1-q^{a}}\sum_{\ell,k=0}^{\infty}\frac{(q^{a};q)_{k+\ell+1}(p^{b};p)_{\ell}}{(p^{c};p)_{\ell}(q^{d};q)_{k+1}(q;q)_{k}(p;p)_{\ell}}x^{k+1}y^{\ell}+\frac{q^{a}}{1-q^{a}}\sum_{\ell,k=0}^{\infty}\frac{(q^{a};q)_{k+\ell}(p^{b};p)_{\ell}}{(p^{c};p)_{\ell}(q^{d};q)_{k}(q;q)_{k}(p;p)_{\ell}}(qx)^{k}y^{\ell}\\
&-\frac{q^{a}}{1-q^{a}}\sum_{\ell,k=0}^{\infty}\frac{(q^{a};q)_{k+\ell}(p^{b};p)_{\ell}}{(p^{c};p)_{\ell}(q^{d};q)_{k}(q;q)_{k}(p;p)_{\ell}}(qx)^{k}(qy)^{\ell}\\
&=\frac{q^{a}x}{1-q^{d}}\Psi_{1}(q^{a+1},p^{b};;p^{c},q^{d+1};q,p,x,y)+\frac{q^{a}}{1-q^{a}}\Psi_{1}(q^{a},p^{b};p^{c},q^{d};q,p,q x,y)\\
&-\frac{q^{a}}{1-q^{a}}\Psi_{1}(q^{a},p^{b};p^{c},q^{d};q,p,q x,q y).
\end{split}
\end{eqnarray*}
The proofs of the relations (\ref{2.2}) and (\ref{2.3}) follow in the same way.
\end{proof}
\begin{thm} The following relations for $\Psi_{1}$ and $\Psi_{2}$ are true
\begin{eqnarray}
\begin{split}
&\Psi_{1}(q^{a},p^{b+1};p^{c},q^{d};q,p,x,y)=\Psi_{1}(q^{a},p^{b};p^{c},q^{d};q,p,x,y)\\
&+\frac{p^{b}(1-q^{a})y}{1-p^{c}}\Psi_{1}(q^{a+1},p^{b+1};p^{c+1},q^{d};q,p,x,y),p^{c}\neq 1,\label{2.4}
\end{split}
\end{eqnarray}
\begin{eqnarray}
\begin{split}
&\Psi_{1}(q^{a},p^{b};p^{c-1},q^{d};q,p,x,y)=\Psi_{1}(q^{a},p^{b};p^{c},q^{d};q,p,x,y)\\
&+\frac{p^{c-1}(1-q^{a})(1-p^{b})y}{(1-p^{c-1})(1-p^{c})}\Psi_{1}(q^{a+1},p^{b+1};p^{c+1},q^{d};q,p,x,y),p^{c},p^{c-1}\neq 1,\\
&\Psi_{1}(q^{a},p^{b};p^{c},q^{d};q,p,x,y)=(1-p^{b})\Psi_{1}(q^{a},p^{b+1};p^{c},q^{d};q,p,x,y)\\
&+p^{b}\Psi_{1}(q^{a},p^{b};p^{c},q^{d};q,p,x,p y)\label{2.5}
\end{split}
\end{eqnarray}
and
\begin{eqnarray}
\begin{split}
&\Psi_{2}(q^{a};p^{b-1},q^{c};q,p,x,y)=\Psi_{2}(q^{a};p^{b},q^{c};q,p,x,y)\\
&+\frac{p^{b-1}(1-q^{a})y}{(1-p^{b-1})(1-p^{b})}\Psi_{2}(q^{a+1};p^{b+1},q^{c};q,p,x,y),p^{b},p^{b-1}\neq 1.\label{2.6}
\end{split}
\end{eqnarray}
\end{thm}
\begin{proof} To prove the relation (\ref{2.4}). Using the relation
\begin{eqnarray*}
\begin{split}
(p^{b};p)_{\ell+1}=(1-p^{b})(p^{b+1};p)_{\ell}=(1-p^{b+\ell})(p^{b};p)_{\ell},
\end{split}
\end{eqnarray*}
 and (\ref{1.5}), we have
\begin{eqnarray*}
\begin{split}
&\Psi_{1}(q^{a},p^{b+1};p^{c},q^{d};q,p,x,y)-\Psi_{1}(q^{a},p^{b};p^{c},q^{d};q,p,x,y)\\
&=p^{b}\sum_{\ell,k=0}^{\infty}\bigg{[}\frac{1-p^{\ell}}{1-p^{b}}\bigg{]}\frac{(q^{a};q)_{k+\ell}(p^{b};p)_{\ell}}{(p^{c};p)_{\ell}(q^{d};q)_{k}(q;q)_{k}(p;p)_{\ell}}x^{k}y^{\ell}\\
&=\frac{p^{b}}{1-p^{b}}\sum_{\ell=1,k=0}^{\infty}\frac{(q^{a};q)_{k+\ell}(p^{b};p)_{\ell}}{(p^{c};p)_{\ell}(q^{d};q)_{k}(q;q)_{k}(p;p)_{\ell-1}}x^{k}y^{\ell}\\
&=\frac{p^{b}(1-q^{a})y}{1-p^{c}}\sum_{\ell,k=0}^{\infty}\frac{(q^{a+1};q)_{k+\ell}(p^{b+1};p)_{\ell}}{(p^{c+1};p)_{\ell}(q^{d};q)_{k}(q;q)_{k}(p;p)_{\ell}}x^{k}y^{\ell+1}\\
&=\frac{p^{b}(1-q^{a})y}{1-p^{c}}\Psi_{1}(q^{a+1},p^{b+1};p^{c+1},q^{d};q,p,x,y),p^{c}\neq1.
\end{split}
\end{eqnarray*}
A similar way to the proof of relation (\ref{2.4}), we obtain the results (\ref{2.5}) and (\ref{2.6})
\end{proof}

\begin{thm} The relations for $\Psi_{1}$ and $\Psi_{2}$ hold true
\begin{eqnarray}
\begin{split}
&(1-q^{a})\Psi_{1}(q^{a+1},p^{b};p^{c},q^{d};q,p,x,\frac{y}{q})+q^{a+1-d}\Psi_{1}(q^{a},p^{b};p^{c},q^{d};q,p,x,y)\\
&=\Psi_{1}(q^{a},p^{b};p^{c},q^{d};q,p,x,\frac{y}{q})+q^{a+1-d}(1-q^{d-1})\Psi_{1}(q^{a},p^{b};q^{d-1};q,p,x,y),\label{2.7}
\end{split}
\end{eqnarray}
\begin{eqnarray}
\begin{split}
&(1-p^{b})\Psi_{1}(q^{a},p^{b+1};p^{c},q^{d};q,p,x,y)=p^{b+1-c}(1-p^{c-1})\Psi_{1}(q^{a},p^{b};p^{c-1},q^{d};q,p,x,y)\\
&+(1-p^{b+1-c})\Psi_{1}(q^{a},p^{b};p^{c},q^{d};q,p,x,y)\label{2.8}
\end{split}
\end{eqnarray}
and
\begin{eqnarray}
\begin{split}
&(1-q^{a})\Psi_{2}(q^{a+1};p^{b},q^{c};q,p,x,\frac{y}{q})+q^{a+1-c}\Psi_{2}(q^{a};p^{b},q^{c};q,p,x,y)\\
&=\Psi_{2}(q^{a};p^{b},q^{c};q,p,x,\frac{y}{q})+q^{a+1-c}(1-q^{c-1})\Psi_{2}(q^{a};p^{b},q^{c-1};q,p,x,y).\label{2.9}
\end{split}
\end{eqnarray}
\end{thm}
\begin{proof}
Using the relationship
\begin{eqnarray*}
\begin{split}
\frac{1-q^{d-1}}{(q^{d-1};q)_{k}}=\frac{1-q^{d+k-1}}{(q^{d};q)_{k}}=\frac{1}{(q^{d};q)_{k-1}},
\end{split}
\end{eqnarray*}
we get
\begin{eqnarray*}
\begin{split}
&q^{a+1-d}(1-q^{d-1})\Psi_{1}(q^{a},p^{b};p^{c},q^{d-1};q,p,x,y)=\sum_{\ell,k=0}^{\infty}\frac{(q^{a+1-d}-q^{a+k})(q^{a};q)_{k+\ell}(p^{b};p)_{\ell}}{(p^{c};p)_{\ell}(q^{d};q)_{k}(q;q)_{k}(p;p)_{\ell}}x^{k}y^{\ell}\\
&=\sum_{\ell,k=0}^{\infty}\frac{(1-q^{a+k+\ell})(q^{a};q)_{k+\ell}(p^{b};p)_{\ell}}{(p^{c};p)_{\ell}(q^{d};q)_{k}(q;q)_{k}(p;p)_{\ell}}x^{k}\bigg{(}\frac{y}{q}\bigg{)}^{\ell}-\sum_{\ell,k=0}^{\infty}\frac{(q^{-\ell}-q^{a+1-d})(q^{a};q)_{k+\ell}(p^{b};p)_{\ell}}{(p^{c};p)_{\ell}(q^{d};q)_{k}(q;q)_{k}(p;p)_{\ell}}x^{k}y^{\ell}\\
&=(1-q^{a})\Psi_{1}(q^{a+1},p^{b};p^{c},q^{d};q,p,x,\frac{y}{q})-\Psi_{1}(q^{a},p^{b};p^{c},q^{d};q,p,x,\frac{y}{q})\\
&+q^{a+1-d}\Psi_{1}(q^{a},p^{b};p^{c},q^{d};q,p,x,y).
\end{split}
\end{eqnarray*}
A similar argument, we obtain the relations (\ref{2.8}) and (\ref{2.9})
\end{proof}
\begin{thm} The bibasic Humbert functions $\Psi_{1}$ and $\Psi_{2}$ satisfy the $q$ and $p$-difference equations
\begin{eqnarray}
\begin{split}
\mathbb{D}_{x,q}^{r}&\Psi_{1}(q^{a},p^{b};p^{c},q^{d};q,p,x,y)=\frac{(q^{a};q)_{r}}{(1-q)^{r}(q^{d};q)_{r}}\Psi_{1}(q^{a+r},p^{b};p^{c},q^{d+r};q,p,x,y),\label{2.10}
\end{split}
\end{eqnarray}
\begin{eqnarray}
\begin{split}
\mathbb{D}_{y,p}^{s}\Psi_{1}(q^{a},p^{b};p^{c},q^{d};q,p,x,y)=&\frac{(q^{a};q)_{s}(p^{b};p)_{r}}{(1-p)^{s}(p^{c};p)_{s}}\Psi_{1}(q^{a+s},p^{b+s};p^{c+s},q^{d};q,p,x,y),\\
\mathbb{D}_{x,p}^{r}\mathbb{D}_{y,p}^{s}\Psi_{1}(q^{a},p^{b};p^{c},q^{d};q,p,x,y)=&\frac{(q^{a};q)_{r+s}(p^{b};p)_{r}}{(1-q)^{r}(1-p)^{s}(p^{c};p)_{s}(q^{d};q)_{r}}\\
&\times \Psi_{1}(q^{a+r+s},p^{b+s};p^{c+s},q^{d+r};q,p,x,y)\label{2.11}
\end{split}
\end{eqnarray}
and
\begin{eqnarray}
\begin{split}
\mathbb{D}_{x,q}^{r}&\Psi_{2}(q^{a};p^{b},q^{c};q,p,x,y)=\frac{(q^{a};q)_{r}}{(1-q)^{r}(q^{c};q)_{r}}\Psi_{2}(q^{a+r};p^{b},q^{c+r};q,p,x,y),\\
\mathbb{D}_{y,p}^{s}&\Psi_{2}(q^{a};p^{b},q^{c};q,p,x,y)=\frac{(q^{a};q)_{s}}{(1-p)^{s}(p^{b};p)_{s}}\Psi_{2}(q^{a+s};p^{b+s},q^{c};q,p,x,y),\\
\mathbb{D}_{x,p}^{r}\mathbb{D}_{y,p}^{s}&\Psi_{2}(q^{a};p^{b},q^{c};q,p,x,y)=\frac{(q^{a};q)_{r+s}\Psi_{2}(q^{a+r+s};p^{b+s},q^{c+r};q,p,x,y)}{(1-q)^{r}(1-p)^{s}(p^{b};p)_{s}(q^{c};q)_{r}}.\label{2.12}
\end{split}
\end{eqnarray}
\end{thm}
\begin{proof}
Calculating the $q$-derivative of both sides of (\ref{1.5}) with respect to the variable $x$, we get
\begin{eqnarray}
\begin{split}
\mathbb{D}_{x,q}&\Psi_{1}(q^{a},p^{b};p^{c},q^{d};q,p,x,y)=\sum_{\ell,k=0}^{\infty}\frac{1-q^{k}}{1-q}\frac{(q^{a};q)_{k+\ell}(p^{b};p)_{\ell}}{(p^{c};p)_{\ell}(q^{d};q)_{k}(q;q)_{k}(p;p)_{\ell}}x^{k-1}y^{\ell}\\
&=\frac{(1-q^{a})}{(1-q)(1-q^{d})}\sum_{\ell,k=0}^{\infty}\frac{(q^{a+1};q)_{k+\ell}(p^{b};p)_{\ell}}{(p^{c};p)_{\ell}(q^{d+1};q)_{k}(q;q)_{k}(p;p)_{\ell}}x^{k}y^{\ell}\\
&=\frac{(1-q^{a})}{(1-q)(1-q^{d})}\Psi_{1}(q^{a+1},p^{b};p^{c},q^{d+1};q,p,x,y).\label{2.13}
\end{split}
\end{eqnarray}
Iterating this $q$-derivative on $\Psi_{1}$ for $r$-times, we obtain (\ref{2.10}).

Using the $p$-derivative in (\ref{1.5}), we get
\begin{eqnarray}
\begin{split}
\mathbb{D}_{y,p}&\Psi_{1}(q^{a},p^{b};p^{c},q^{d};q,p,x,y)=\sum_{\ell=1,k=0}^{\infty}\frac{1}{1-p}\frac{(q^{a};q)_{k+\ell}(p^{b};p)_{\ell}}{(p^{c};p)_{\ell}(q^{d};q)_{k}(p;p)_{\ell-1}(q;q)_{k}}x^{k}y^{\ell-1}\\
&=\frac{(1-q^{a})(1-p^{b})}{(1-p)(1-p^{c})}\sum_{\ell,k=0}^{\infty}\frac{(q^{a+1};q)_{k+\ell}(p^{b+1};p)_{\ell}}{(p^{c+1};p)_{\ell}(q^{d};q)_{k}(q;q)_{k}(p;p)_{\ell}}x^{k}y^{\ell}\\
&=\frac{(1-q^{a})(1-p^{b})}{(1-p)(1-p^{c})}\Psi_{1}(q^{a+1},p^{b+1};p^{c+1};q,p,x,y).\label{2.14}
\end{split}
\end{eqnarray}
Iterating this $p$-derivative on $\Psi_{1}$ for $s$-times, we obtain (\ref{2.11}). Similarly, the $q$-derivatives given by (\ref{1.4}) can be proved (\ref{2.12}).
\end{proof}
\begin{thm} The $q$-differential relations for $\Psi_{1}$ and $\Psi_{2}$ hold
\begin{eqnarray}
\begin{split}
x\mathbb{D}_{x,q}\Psi_{1}(q^{a},p^{b};p^{c},q^{d};q,p,x,y)=&\frac{(1-q^{d-1})}{(1-q)q^{d-1}}\bigg{[}\Psi_{1}(q^{a},p^{b};p^{c},q^{d-1};q,p,x,y)\\
&-\Psi_{1}(q^{a},p^{b};p^{c},q^{d};q,p,x,y)\bigg{]},\label{2.15}
\end{split}
\end{eqnarray}
\begin{eqnarray}
\begin{split}
y\mathbb{D}_{y,p}\Psi_{1}(q^{a},p^{b};p^{c},q^{d};q,p,x,y)=&\frac{1-p^{b}}{(1-p)p^{b}}\bigg{[}\Psi_{1}(q^{a},p^{b+1};p^{c},q^{d};q,p,x,y)\\
&-\Psi_{1}(q^{a},p^{b};p^{c},q^{d};q,p,x,y)\bigg{]},\\
y\mathbb{D}_{y,p}\Psi_{1}(q^{a},p^{b};p^{c},q^{d};q,p,x,y)=&\frac{1-p^{c-1}}{(1-p)p^{c-1}}\bigg{[}\Psi_{1}(q^{a},p^{b};p^{c-1},q^{d};q,p,x,y)\\
&-\Psi_{1}(q^{a},p^{b};p^{c},q^{d};q,p,x,y)\bigg{]}\label{2.16}
\end{split}
\end{eqnarray}
and
\begin{eqnarray}
\begin{split}
x\mathbb{D}_{x,q}\Psi_{2}(q^{a};p^{b},q^{c};q,p,x,y)=&\frac{(1-q^{c-1})}{(1-q)q^{c-1}}\bigg{[}\Psi_{2}(q^{a};p^{b},q^{c-1};q,p,x,y)\\
&-\Psi_{2}(q^{a};p^{b},q^{c};q,p,x,y)\bigg{]},\\
y\mathbb{D}_{y,p}\Psi_{2}(q^{a};p^{b},q^{c};q,p,x,y)=&\frac{1-p^{b-1}}{(1-p)p^{b-1}}\bigg{[}\Psi_{2}(q^{a};p^{b-1},q^{c};q,p,x,y)\\
&-\Psi_{2}(q^{a};p^{b},q^{c};q,p,x,y)\bigg{]}.\label{2.17}
\end{split}
\end{eqnarray}
\end{thm}
\begin{proof}
Multiplying (\ref{2.13}) by $x$ and substituting the value of $\Psi_{1}(q^{a+1},p^{b};p^{c},q^{d+1};q,p,x,y)$ from (\ref{2.2}) into (\ref{2.13}), we get (\ref{2.15}).

The proofs (\ref{2.16})-(\ref{2.17}) are similar to the proof of the corresponding identity (\ref{2.15}) for the bibasic Humbert functions $\Psi_{1}$ and $\Psi_{2}$ are omitted most steps and give only outlines wherever necessary.
\end{proof}
\begin{thm}
The following relations for $\Psi_{1}$ and $\Psi_{2}$ hold:
\begin{eqnarray}
\begin{split}
&(1-q^{d-1})\Psi_{1}(q^{a},p^{b};p^{c},q^{d-1};q,p,x,y)=(1-q)x\mathbb{D}_{x,q}\Psi_{1}(q^{a},p^{b};p^{c},q^{d};q,p,x,y)\\
&+(1-q^{d-1})\Psi_{1}(q^{a},p^{b};p^{c},q^{d};q,p,q x,y),\label{2.18}
\end{split}
\end{eqnarray}
\begin{eqnarray}
\begin{split}
&(1-p^{c-1})\Psi_{1}(q^{a},p^{b};p^{c-1},q^{d};q,p,x,y)=(1-p)y\mathbb{D}_{y,p}\Psi_{1}(q^{a},p^{b};p^{c},q^{d};q,p,x,y)\\
&+(1-p^{c-1})\Psi_{1}(q^{a},p^{b};p^{c},q^{d};q,p,x,p y),\\
&(1-p^{b})\Psi_{1}(q^{a},p^{b+1};p^{c},q^{d};q,p,x,y)=(1-p)y\mathbb{D}_{y,p}\Psi_{1}(q^{a},p^{b};p^{c},q^{d};q,p,x,y)\\
&+(1-p^{b})\Psi_{1}(q^{a},p^{b};p^{c},q^{d};q,p,x,p y),\label{2.19}
\end{split}
\end{eqnarray}
\begin{eqnarray}
\begin{split}
&(1-q^{a})\Psi_{1}(q^{a+1},p^{b};p^{c},q^{d};q,p,x,xy)=(1-q^{a})\Psi_{1}(q^{a},p^{b};p^{c},q^{d};q,p,x,xy)\\
&+(1-q)q^{a}x\mathbb{D}_{x,q}\Psi_{1}(q^{a},p^{b};p^{c},q^{d};q,p,x,xy),\\
&(1-q^{a})\Psi_{1}(q^{a+1},p^{b};p^{c},q^{d};q,p,x,xy)=(1-q)x\mathbb{D}_{x,q}\Psi_{1}(q^{a},p^{b};p^{c},q^{d};q,p,x,xy)\\
&+(1-q^{a})\Psi_{1}(q^{a},p^{b};p^{c},q^{d};q,p,q x,q xy),\label{2.20}
\end{split}
\end{eqnarray}
\begin{eqnarray}
\begin{split}
&(1-q^{a})\Psi_{1}(q^{a+1},p^{b};p^{c},q^{d};q,p,xy,y)=(1-q^{a})\Psi_{1}(q^{a},p^{b};p^{c},q^{d};q,p,xy,y)\\
&+(1-q)q^{a}y\mathbb{D}_{y,q}\Psi_{1}(q^{a},p^{b};p^{c},q^{d};q,p,xy,y),\\
&(1-q^{a})\Psi_{1}(q^{a+1},p^{b};p^{c},q^{d};q,p,xy,y)=(1-q)y\mathbb{D}_{y,q}\Psi_{1}(q^{a},p^{b};p^{c},q^{d};q,p,xy,y)\\
&+(1-q^{a})\Psi_{1}(q^{a},p^{b};p^{c},q^{d};q,p,q xy,q y),\label{2.21}
\end{split}
\end{eqnarray}
\begin{eqnarray}
\begin{split}
&(1-q^{c-1})\Psi_{2}(q^{a};p^{c},q^{c-1};q,p,x,y)=(1-q)x\mathbb{D}_{x,q}\Psi_{2}(q^{a};p^{b},q^{c};q,p,x,y)\\
&+(1-q^{c-1})\Psi_{2}(q^{a};p^{b},q^{c};q,p,q x,y),\\
&(1-p^{b-1})\Psi_{2}(q^{a};p^{b-1},q^{c};q,p,x,y)=(1-p)y\mathbb{D}_{y,p}\Psi_{2}(q^{a};p^{b},q^{c};q,p,x,y)\\
&+(1-p^{b-1})\Psi_{2}(q^{a};p^{b},q^{c};q,p,x,p y),\label{2.22}
\end{split}
\end{eqnarray}
\begin{eqnarray}
\begin{split}
&(1-q^{a})\Psi_{2}(q^{a+1};p^{b},q^{c};q,p,x,xy)=(1-q^{a})\Psi_{2}(q^{a};p^{b},q^{c};q,p,x,xy)\\
&+(1-q)q^{a}x\mathbb{D}_{x,q}\Psi_{2}(q^{a};p^{b},q^{c};q,p,x,xy),\\
&(1-q^{a})\Psi_{2}(q^{a+1};p^{b},q^{c};q,p,x,xy)=(1-q)x\mathbb{D}_{x,q}\Psi_{2}(q^{a};p^{b},q^{c};q,p,x,xy)\\
&+(1-q^{a})\Psi_{2}(q^{a};p^{b},q^{c};q,p,q x,q xy)\label{2.23}
\end{split}
\end{eqnarray}
and
\begin{eqnarray}
\begin{split}
&(1-q^{a})\Psi_{2}(q^{a+1};p^{b},q^{c};q,p,xy,y)=(1-q^{a})\Psi_{2}(q^{a};p^{b},q^{c};q,p,xy,y)\\
&+(1-q)q^{a}y\mathbb{D}_{y,q}\Psi_{2}(q^{a};p^{b},q^{c};q,p,xy,y),\\
&(1-q^{a})\Psi_{2}(q^{a+1};p^{b},q^{c};q,p,xy,y)=(1-q)y\mathbb{D}_{y,q}\Psi_{2}(q^{a};p^{b},q^{c};q,p,xy,y)\\
&+(1-q^{a})\Psi_{2}(q^{a};p^{b},q^{c};q,p,q xy,q y).\label{2.24}
\end{split}
\end{eqnarray}
\end{thm}
\begin{proof}
Using (\ref{1.5}) and (\ref{2.13}), we obtain
\begin{eqnarray*}
\begin{split}
&(1-q^{d-1})\Psi_{1}(q^{a},p^{b};p^{c},q^{d-1};q,p,x,y)=\sum_{\ell,k=0}^{\infty}\frac{(q^{a};q)_{k+\ell}(p^{b};p)_{\ell}}{(p^{c};p)_{\ell}(q^{d};q)_{k-1}(q;q)_{k}(p;p)_{\ell}}x^{k}y^{\ell}\\
&=\sum_{\ell,k=0}^{\infty}\frac{(1-q^{d+k-1})(q^{a};q)_{k+\ell}(p^{b};p)_{\ell}}{(p^{c};p)_{\ell}(q^{d};q)_{k}(q;q)_{k}(p;p)_{\ell}}x^{k}y^{\ell}\\
&=\sum_{\ell,k=0}^{\infty}\frac{(1-q^{k}+q^{k}(1-q^{d-1}))(q^{a};q)_{k+\ell}(p^{b};p)_{\ell}}{(p^{c};p)_{\ell}(q^{d};q)_{k}(q;q)_{k}(p;p)_{\ell}}x^{k}y^{\ell}\\
&=(1-q)x\mathbb{D}_{x,q}\Psi_{1}(q^{a},p^{b};p^{c},q^{d};q,p,x,y)+(1-q^{d-1})\Psi_{1}(q^{a},p^{b};p^{c},q^{d};q,p,q x,y).
\end{split}
\end{eqnarray*}
The proofs results (\ref{2.19})-(\ref{2.24}) for $\Psi_{1}$ and $\Psi_{2}$ follows similarly from the identity (\ref{2.18}), and is omitted most steps and give only outlines wherever necessary.
\end{proof}
\begin{thm}
The $q$-derivatives of bibasic Humbert functions $\Psi_{1}$ and $\Psi_{2}$ with respect to their parameters satisfies the relations
\begin{eqnarray}
\begin{split}
\mathbb{D}_{a,q}\Psi_{1}(q^{a},p^{b};p^{c},q^{d};q,p,x,y)=&-\frac{1}{1-q^{a}}\bigg{[}x\mathbb{D}_{x,q}\Psi_{1}+\frac{1-p}{1-q}y\mathbb{D}_{y,p}\\
&+\frac{1}{1-q}\Psi_{1}(q x,p y)-\frac{1}{1-q}\Psi_{1}(q x,q y)\bigg{]},q^{a}\neq1,\label{2.25}
\end{split}
\end{eqnarray}
\begin{eqnarray}
\begin{split}
\mathbb{D}_{a,q}\Psi_{1}(q^{a},p^{b};p^{c},q^{d};q,p,x,y)=&-\frac{1}{1-q^{a}}\bigg{[}\frac{1-p}{1-q}y\mathbb{D}_{y,p}\Psi_{1}+x\mathbb{D}_{x,q}\Psi_{1}(q y)\\
&+\frac{1}{1-q}\Psi_{1}(p y)-\frac{1}{1-q}\Psi_{1}(q y)\bigg{]},q^{a}\neq1,\\
\mathbb{D}_{d,q}\Psi_{1}(q^{a},p^{b};p^{c},q^{d};q,p,x,y)=&\frac{1}{1-q^{d}}x\mathbb{D}_{x,q}\Psi_{1}(q^{d+1}),q^{d}\neq1,\\
\mathbb{D}_{b,p}\Psi_{1}(q^{a},p^{b};p^{c},q^{d};q,p,x,y)=&-\frac{1}{1-p^{b}}y\mathbb{D}_{y,p}\Psi_{1},p^{b}\neq1,\\
\mathbb{D}_{c,p}\Psi_{1}(q^{a},p^{b};p^{c},q^{d};q,p,x,y)=&\frac{1}{1-p^{c}}y\mathbb{D}_{y,p}\Psi_{1}(p^{c+1}),p^{c}\neq1\label{2.26}
\end{split}
\end{eqnarray}
and
\begin{eqnarray}
\begin{split}
\mathbb{D}_{a,q}\Psi_{2}(q^{a};p^{b},q^{c};q,p,x,y)=&-\frac{1}{1-q^{a}}\bigg{[}x\mathbb{D}_{x,q}\Psi_{2}+\frac{1-p}{1-q}y\mathbb{D}_{y,p}\Psi_{2}(q x)\\
&+\frac{1}{1-q}\Psi_{2}(q x,p y)-\frac{1}{1-q}\Psi_{2}(q x,q y)\bigg{]},q^{a}\neq1,\\
\mathbb{D}_{a,q}\Psi_{2}(q^{a};p^{b},q^{c};q,p,x,y)=&-\frac{1}{1-q^{a}}\bigg{[}\frac{1-p}{1-q}y\mathbb{D}_{y,p}\Psi_{2}+x\mathbb{D}_{x,q}\Psi_{2}(q y)\\
&+\frac{1}{1-q}\Psi_{2}(p y)-\frac{1}{1-q}\Psi_{2}(q y)\bigg{]},q^{a}\neq1,\\
\mathbb{D}_{c,q}\Psi_{2}(q^{a};p^{b},q^{c};q,p,x,y)=&\frac{1}{1-q^{c}}x\mathbb{D}_{x,q}\Psi_{2}(q^{c+1}),q^{c}\neq1,\\
\mathbb{D}_{b,p}\Psi_{2}(q^{a};p^{b},q^{c};q,p,x,y)=&\frac{1}{1-p^{b}}y\mathbb{D}_{y,p}\Psi_{2}(p^{b+1}),p^{b}\neq1.\label{2.27}
\end{split}
\end{eqnarray}
\end{thm}
\begin{proof}
Calculating the $q$-derivative of both sides of (\ref{1.5}) with respect to the variable $a$, we get
\begin{eqnarray*}
\begin{split}
\mathbb{D}_{a,q}&\Psi_{1}(q^{a},p^{b};p^{c},q^{d};q,p,x,y)=\sum_{\ell,k=0}^{\infty}\frac{(q^{a};q)_{k+\ell}-(q^{a+1};q)_{k+\ell}}{(1-q)q^{a}}\frac{(p^{b};p)_{\ell}}{(p^{c};p)_{\ell}(q^{d};q)_{k}(q;q)_{k}(p;p)_{\ell}}x^{k}y^{\ell}\\
&=-\frac{1}{1-q^{a}}\sum_{\ell,k=0}^{\infty}\frac{1-q^{k+\ell}}{1-q}\frac{(q^{a};q)_{k+\ell}(p^{b};p)_{\ell}}{(p^{c};p)_{\ell}(q^{d};q)_{k}(q;q)_{k}(p;p)_{\ell}}x^{k}y^{\ell}\\
&=-\frac{1}{1-q^{a}}\sum_{\ell,k=0}^{\infty}\bigg{[}\frac{1-q^{k}}{1-q}+\frac{q^{k}(1-p)}{1-q}\frac{1-p^{\ell}}{1-p}+q^{k}\frac{p^{\ell}-q^{\ell}}{1-q}\bigg{]}\frac{(q^{a};q)_{k+\ell}(p^{b};p)_{\ell}}{(p^{c};p)_{\ell}(q^{d};q)_{k}(q;q)_{k}(p;p)_{\ell}}x^{k}y^{\ell}\\
&=-\frac{1}{1-q^{a}}\bigg{[}x\mathbb{D}_{x,q}\Psi_{1}+\frac{1-p}{1-q}y\mathbb{D}_{y,p}\Psi_{1}(q x)+\frac{1}{1-q}\Psi_{1}(q x,p y)-\frac{1}{1-q}\Psi_{1}(q x,q y)\bigg{]}.
\end{split}
\end{eqnarray*}
The proofs (\ref{2.26})-(\ref{2.27}) for $\Psi_{1}$ and $\Psi_{2}$ are similar to the proof of results (\ref{2.25}) and are omitted.
\end{proof}
\begin{thm}
For $r\in \mathbb{N}$, the differentiation formulas for $\Psi_{1}$ and $\Psi_{2}$ hold true
\begin{eqnarray}
\begin{split}
\mathbb{D}_{y,p}^{r}\bigg{[}y^{b+r-1}\Psi_{1}(q^{a},p^{b};p^{c},q^{d};q,p,x,y)\bigg{]}=\frac{(p^{b};p)_{r}}{(1-p)^{r}}y^{b-1}\Psi_{1}(q^{a},p^{b+r};p^{c},q^{d};q,p,x,y),\label{2.28}
\end{split}
\end{eqnarray}
\begin{eqnarray}
\begin{split}
\mathbb{D}_{x,q}^{r}\bigg{[}x^{a+r-1}\Psi_{1}(q^{a},p^{b};p^{c},q^{d};q,p,x,xy)\bigg{]}=\frac{(q^{a};q)_{r}}{(1-q)^{r}}x^{a-1}\Psi_{1}(q^{a+r},p^{b};p^{c},q^{d};q,p,x,xy),\\
\mathbb{D}_{y,q}^{r}\bigg{[}y^{a+r-1}\Psi_{1}(q^{a},p^{b};p^{c},q^{d};q,p,xy,y)\bigg{]}=\frac{(q^{a};q)_{r}}{(1-q)^{r}}y^{a-1}\Psi_{1}(q^{a+r},p^{b};p^{c},q^{d};q,p,xy,y)\label{2.29}
\end{split}
\end{eqnarray}
and
\begin{eqnarray}
\begin{split}
\mathbb{D}_{x,q}^{r}\bigg{[}x^{a+r-1}\Psi_{2}(q^{a};p^{b},q^{c};q,p,x,xy)\bigg{]}=\frac{(q^{a};q)_{r}}{(1-q)^{r}}x^{a-1}\Psi_{2}(q^{a+r};p^{b},q^{c};q,p,x,xy),\\
\mathbb{D}_{y,q}^{r}\bigg{[}y^{a+r-1}\Psi_{2}(q^{a};p^{b},q^{c};q,p,xy,y)\bigg{]}=\frac{(q^{a};q)_{r}}{(1-q)^{r}}y^{a-1}\Psi_{2}(q^{a+r};p^{b},q^{c};q,p,xy,y).\label{2.30}
\end{split}
\end{eqnarray}
\end{thm}
\begin{proof}
Using
\begin{eqnarray*}
\begin{split}
&\mathbb{D}_{y,p}^{r}\bigg{[}y^{b+\ell+r-1}\bigg{]}=\frac{(p^{b+\ell};p)_{r}}{(1-p)^{r}}y^{b+\ell-1}
\end{split}
\end{eqnarray*}
and
\begin{eqnarray*}
\begin{split}
(p^{b};p)_{\ell}(p^{b+\ell};p)_{r}=(p^{b};p)_{\ell+r}=(p^{b};p)_{r}(p^{b+r};p)_{\ell},
\end{split}
\end{eqnarray*}
we get 
\begin{eqnarray*}
\begin{split}
&\mathbb{D}_{y,p}^{r}\bigg{[}y^{b+r-1}\Psi_{1}(q^{a},p^{b};p^{c},q^{d};q,p,x,y)\bigg{]}=\sum_{\ell,k=0}^{\infty}\frac{(q^{a};q)_{k+\ell}(p^{b};p)_{\ell}}{(p^{c};p)_{\ell}(q^{d};q)_{k}(q;q)_{k}(p;p)_{\ell}}x^{k}\mathbb{D}_{y,p}^{r}\bigg{[}y^{b+\ell+r-1}\bigg{]}\\
&=\frac{1}{(1-p)^{r}}y^{b-1}\sum_{\ell,k=0}^{\infty}\frac{(q^{a};q)_{k+\ell}(p^{b};p)_{\ell}(p^{b+\ell};p)_{r}}{(p^{c};p)_{\ell}(q^{d};q)_{k}(q;q)_{k}(p;p)_{\ell}}x^{k}y^{\ell}\\
&=\frac{(p^{b};p)_{r}}{(1-p)^{r}}y^{b-1}\sum_{\ell,k=0}^{\infty}\frac{(q^{a};q)_{k+\ell}(p^{b+r};p)_{\ell}}{(p^{c};p)_{\ell}(q^{d};q)_{k}(q;q)_{k}(p;p)_{\ell}}x^{k}y^{\ell}\\
&=\frac{(p^{b};p)_{r}}{(1-p)^{r}}y^{b-1}\Psi_{1}(q^{a},p^{b+r};p^{c},q^{d};q,p,x,y).
\end{split}
\end{eqnarray*}
The proof Eqs (\ref{2.29})-(\ref{2.30}) are on the same lines as of Eq. (\ref{2.28}) .
\end{proof}
\begin{thm} The summation formulas for $\Psi_{1}$ and $\Psi_{2}$ hold true
\begin{eqnarray}
\begin{split}
\Psi_{1}(q^{a},p^{b};p^{c},q^{d};q,p,x,y)=\sum_{\ell=0}^{\infty}\frac{(q^{a};q)_{\ell}(p^{b};p)_{\ell}}{(p^{c};p)_{\ell}(p;p)_{\ell}}y^{\ell}\;_{2}\phi_{1}(q^{a+\ell},0;q^{d};q,x)\label{2.31}
\end{split}
\end{eqnarray}
and
\begin{eqnarray}
\begin{split}
\Psi_{2}(q^{a};p^{b},q^{c};q,p,x,y)=\sum_{\ell=0}^{\infty}\frac{(q^{a};q)_{\ell}}{(p^{b};p)_{\ell}(p;p)_{\ell}}y^{\ell}\;_{2}\phi_{1}(q^{a+\ell},0;q^{c};q,x).\label{2.32}
\end{split}
\end{eqnarray}
\end{thm}
\begin{proof}
For any two integers $\ell$ and $k$, the following holds
\begin{eqnarray*}
\begin{split}
(q^{a};q)_{\ell+k}=(q^{a};q)_{\ell}(q^{a+\ell};q)_{k}=(q^{a};q)_{k}(q^{a+k};q)_{\ell},
\end{split}
\end{eqnarray*}
we have
\begin{eqnarray*}
\begin{split}
&\Psi_{1}(q^{a},p^{b};p^{c},q^{d};q,p,x,y)=\sum_{\ell,k=0}^{\infty}\frac{(q^{a};q)_{\ell}(q^{a+\ell};q)_{k}(p^{b};p)_{\ell}}{(p^{c};p)_{\ell}(q^{d};q)_{k}(q;q)_{k}(p;p)_{\ell}}x^{k}y^{\ell}\\
&=\sum_{\ell=0}^{\infty}\frac{(q^{a};q)_{\ell}(p^{b};p)_{\ell}}{(p^{c};p)_{\ell}(p;p)_{\ell}}y^{\ell}\sum_{k=0}^{\infty}\frac{(q^{a+\ell};q)_{k}}{(q^{d};q)_{k}(q;q)_{k}}x^{k}=\sum_{\ell=0}^{\infty}\frac{(q^{a};q)_{\ell}(p^{b};p)_{\ell}}{(p^{c};p)_{\ell}(p;p)_{\ell}}y^{\ell}\;_{2}\phi_{1}(q^{a+\ell},0;q^{d};q,x),
\end{split}
\end{eqnarray*}
we obtain (\ref{2.32}) can be proved similarly.
\end{proof}
\begin{thm} The connection relation between bibasic Humbert functions $\Psi_{1}$ and $\Psi_{2}$ is shown as follows
\begin{eqnarray}
\begin{split}
&\Psi_{1}(q^{a},p^{b};p^{c},q^{d};q,p,x,y)=\frac{(p^{b};p)_{\infty}}{(p^{c};p)_{\infty}}\sum_{s=0}^{\infty}\frac{(p^{c-b};p)_{s}p^{bs}}{(p;p)_{s}}\Psi_{2}(q^{a};0,q^{d};q,p,x,p^{s}y).\label{2.33}
\end{split}
\end{eqnarray}
\end{thm}
\begin{proof}
Using the identities
\begin{eqnarray*}
\begin{split}
(p^{b};p)_{\ell}=\frac{(p^{b};p)_{\infty}}{(p^{b+\ell};q)_{\infty}},\\
(p^{c};p)_{\ell}=\frac{(p^{c};p)_{\infty}}{(p^{c+\ell};q)_{\infty}},
\end{split}
\end{eqnarray*}
and
\begin{eqnarray*}
\begin{split}
\frac{(p^{c+\ell};p)_{\infty}}{(p^{b+\ell};p)_{\infty}}=\sum_{s=0}^{\infty}\frac{(p^{c-b};p)_{s}}{(p;p)_{s}}\bigg{(}p^{b+\ell}\bigg{)}^{s},
\end{split}
\end{eqnarray*}
substitution of this gives 
\begin{eqnarray*}
\begin{split}
&\Psi_{1}(q^{a},p^{b};p^{c},q^{d};q,p,x,y)=\frac{(p^{b};p)_{\infty}}{(p^{c};p)_{\infty}}\sum_{\ell,k=0}^{\infty}\frac{(q^{a};q)_{k+\ell}(p^{c+\ell};p)_{\infty}}{(p^{b+\ell};p)_{\infty}(q^{d};q)_{k}(q;q)_{k}(p;p)_{\ell}}x^{k}y^{\ell}\\
&=\frac{(p^{b};p)_{\infty}}{(p^{c};p)_{\infty}}\sum_{\ell,k,s=0}^{\infty}\frac{(q^{a};q)_{k+\ell}(p^{c-b};p)_{s}}{(q^{d};q)_{k}(q;q)_{k}(p;p)_{\ell}(p;p)_{s}}\bigg{(}p^{b+\ell}\bigg{)}^{s}x^{k}y^{\ell}\\
&=\frac{(p^{b};p)_{\infty}}{(p^{c};p)_{\infty}}\sum_{s=0}^{\infty}\frac{(p^{c-b};p)_{s}p^{bs}}{(p;p)_{s}}\sum_{\ell,k=0}^{\infty}\frac{(q^{a};q)_{k+\ell}}{(q^{d};q)_{k}(q;q)_{k}(p;p)_{\ell}}p^{s\ell}x^{k}y^{\ell}\\
&=\frac{(p^{b};p)_{\infty}}{(p^{c};p)_{\infty}}\sum_{s=0}^{\infty}\frac{(p^{c-b};p)_{s}p^{bs}}{(p;p)_{s}}\Psi_{2}(q^{a};0,q^{d};q,p,x,p^{s}y).
\end{split}
\end{eqnarray*}
\end{proof}
\begin{thm}
The series representations for $\Psi_{1}$ and $\Psi_{2}$ in (\ref{1.5})- (\ref{1.6}) satisfy the results
\begin{eqnarray}
\begin{split}
\Psi_{1}(q^{a},p^{b};p^{c},q^{d};q,p,x,y)=&\frac{(q^{a};q)_{\infty}(p^{b};p)_{\infty}}{(q^{d};q)_{\infty}(p^{c};p)_{\infty}}\sum_{r,s,k,\ell=0}^{\infty}\frac{(q^{a+k};q)_{\ell}(q^{d-a};q)_{r}(p^{c-b};p)_{s}}{(q;q)_{r}(p;p)_{s}(q;q)_{k}(p;p)_{\ell}}\\
&\times q^{ar}p^{bs}\big{(}q^{r}x\big{)}^{k}\big{(}p^{s}y\big{)}^{\ell},\label{2.34}
\end{split}
\end{eqnarray}

\begin{eqnarray}
\begin{split}
\Psi_{1}(q^{a},p^{b};p^{c},q^{d};q,p,x,y)&=\frac{(q^{a};q)_{\infty}(p^{b};p)_{\infty}(xq^{a};q)_{\infty}}{(q^{d};q)_{\infty}(p^{c};p)_{\infty}(x;q)_{\infty}}\\
&\times \sum_{r,s,\ell=0}^{\infty}\frac{(q^{d-a};q)_{r}(p^{c-b};p)_{s}q^{ar}p^{bs}}{(q;q)_{r}(p;p)_{s}(xq^{a};q)_{\ell}(p;p)_{\ell}}\big{(}q^{r}p^{s}y\big{)}^{\ell},\label{2.35}
\end{split}
\end{eqnarray}
\begin{eqnarray}
\begin{split}
\Psi_{1}(q^{a},p^{b};p^{c},q^{d};q,p,x,y)&=\frac{(q^{a};q)_{\infty}(p^{b};p)_{\infty}}{(q^{d};q)_{\infty}(p^{c};p)_{\infty}}\sum_{r,s,\ell=0}^{\infty}\frac{(q^{d-a};q)_{r}(p^{c-b};p)_{s}q^{ar}p^{bs}}{(q;q)_{r}(p;p)_{s}(p;p)_{\ell}}\\
&\times \big{(}p^{s}y\big{)}^{\ell}\:_{1}\mathbf{\Phi}_{0}(q^{a+\ell};-;q,q^{\ell}x)\label{2.36}
\end{split}
\end{eqnarray}
and
\begin{eqnarray}
\begin{split}
&\Psi_{2}(q^{a};p^{b},q^{c};q,p,x,y)=\frac{(q^{a};q)_{\infty}}{(q^{c};q)_{\infty}}\sum_{r,\ell,k=0}^{\infty}\frac{(q^{a+k};q)_{\ell}(q^{c-a};q)_{r}}{(p^{b};p)_{\ell}(q;q)_{r}(q;q)_{k}(p;p)_{\ell}}q^{(a+k)r}x^{k}y^{\ell}.\label{2.37}
\end{split}
\end{eqnarray}
\end{thm}
\begin{proof}
We can write the series of the function $\Psi_{1}$ as
\begin{eqnarray*}
\begin{split}
&\Psi_{1}(q^{a},p^{b};p^{c},q^{d};q,p,x,y)=\sum_{\ell,k=0}^{\infty}\frac{(q^{a};q)_{k}(q^{a+k};q)_{\ell}(p^{b};p)_{\ell}}{(p^{c};p)_{\ell}(q^{d};q)_{k}(q;q)_{k}(p;p)_{\ell}}x^{k}y^{\ell}\\
&=\sum_{\ell,k=0}^{\infty}\frac{(q^{a};q)_{k}(q^{a+k};q)_{\ell}(p^{b};p)_{\infty}(p^{c+\ell};p)_{\infty}}{(p^{b+\ell};p)_{\infty}(p^{c};p)_{\infty}(q^{d};q)_{k}(q;q)_{k}(p;p)_{\ell}}x^{k}y^{\ell}\\
&=\frac{(p^{b};p)_{\infty}}{(p^{c};p)_{\infty}}\sum_{\ell,k=0}^{\infty}\frac{(q^{a};q)_{k}(q^{a+k};q)_{\ell}(p^{c+\ell};p)_{\infty}}{(p^{b+\ell};p)_{\infty}(q^{d};q)_{k}(q;q)_{k}(p;p)_{\ell}}x^{k}y^{\ell}\\
&=\frac{(q^{a};q)_{\infty}(p^{b};p)_{\infty}}{(q^{d};q)_{\infty}(p^{c};p)_{\infty}}\sum_{\ell,k=0}^{\infty}\frac{(q^{d+k};q)_{\infty}(q^{a+k};q)_{\ell}(p^{c+\ell};p)_{\infty}}{(p^{b+\ell};p)_{\infty}(q^{a+k};q)_{\infty}(q;q)_{k}(p;p)_{\ell}}x^{k}y^{\ell}\\
&=\frac{(q^{a};q)_{\infty}(p^{b};p)_{\infty}}{(q^{d};q)_{\infty}(p^{c};p)_{\infty}}\sum_{\ell,k,s,r=0}^{\infty}\frac{(q^{a+k};q)_{\ell}(q^{d-a};q)_{r}(p^{c-b};p)_{s}}{(q;q)_{k}(p;p)_{\ell}(q;q)_{r}(p;p)_{s}}\bigg{(}q^{a+k}\bigg{)}^{r}\bigg{(}p^{b+\ell}\bigg{)}^{s}x^{k}y^{\ell}.
\end{split}
\end{eqnarray*}
The formulas (\ref{2.35})-(\ref{2.37}) are proved in a similar manner.
\end{proof}
\begin{thm}
The $q$-integral representations for $\Psi_{1}$ are true
\begin{eqnarray}
\begin{split}
\Psi_{1}(q^{a},p^{b};p^{c},q^{d};q,p,x,y)&=\frac{\Gamma_{p}(c)}{\Gamma_{p}(b)\Gamma_{p}(c-b)}\int_{0}^{1}t^{b-1}\frac{(pt;p)_{\infty}}{(tp^{c-b};p)_{\infty}}\\
&\times \Psi_{1}(q^{a},p^{e};p^{e},q^{d};q,p,x,yt)d_{p}t,\label{2.38}
\end{split}
\end{eqnarray}
\begin{eqnarray}
\begin{split}
\Psi_{1}(q^{a},p^{b};p^{c},q^{d};q,p,x,y)&=\frac{\Gamma_{p}(c)}{\Gamma_{p}(b)\Gamma_{p}(c-b)}\int_{0}^{1}t^{b-1}\frac{(pt;p)_{\infty}}{(tp^{c-b};p)_{\infty}}\\
&\times\Psi_{2}(q^{a};0,q^{d};q,p,x,yt)d_{p}t,\label{2.39}
\end{split}
\end{eqnarray}
\begin{eqnarray}
\begin{split}
\Psi_{1}(q^{a},p^{b};p^{c},q^{d};q,p,x,y)&=\frac{1}{\Gamma_{p}(b)}\int_{0}^{\frac{1}{1-p}}\;E_{p}(-pt)t^{b-1}\\
&\times\Psi_{2}(q^{a};p^{c},q^{d};q,p,x,(1-p)yt)d_{p}t\label{2.40}
\end{split}
\end{eqnarray}
and
\begin{eqnarray}
\begin{split}
\Psi_{1}(q^{a},p^{b};p^{c},q^{d};q,p,x,y)&=\frac{1}{\Gamma_{q}(a)}\int_{0}^{\frac{1}{1-q}}\;E_{q}(-qt)t^{a-1}\:_{1}\Phi_{1}(p^{b};p^{c};(1-q)yt)\\
&\times \:_{0}\mathbf{\Phi}_{1}(-;q^{d};(1-q)xt)d_{q}t.\label{2.41}
\end{split}
\end{eqnarray}
\end{thm}
\begin{proof}
Using
\begin{eqnarray*}
\begin{split}
\frac{(p^{b};p)_{\ell}}{(p^{c};p)_{\ell}}=\frac{\Gamma_{p}(c)}{\Gamma_{p}(b)\Gamma_{p}(c-b)}\int_{0}^{1}t^{b+\ell-1}\frac{(pt;p)_{\infty}}{(tp^{c-b};p)_{\infty}}d_{p}t,
\end{split}
\end{eqnarray*}
for $0<\Re(b)<\Re(c)$, $c-b\neq 0, -1, -2, -3, \ldots$, $Re(b)>0$ and $\ell\geq 0$, we obtain an $p$-integral representation for $\Psi_{1}$
\begin{eqnarray*}
\begin{split}
&\Psi_{1}(q^{a},p^{b};p^{c},q^{d};q,p,x,y)=\sum_{\ell,k=0}^{\infty}\frac{(q^{a};q)_{k+\ell}(p^{b};p)_{\ell}}{(p^{c};p)_{\ell}(q^{d};q)_{k}(q;q)_{k}(p;p)_{\ell}}x^{k}y^{\ell}\\
&=\frac{\Gamma_{p}(c)}{\Gamma_{p}(b)\Gamma_{p}(c-b)}\sum_{\ell,k=0}^{\infty}\frac{(q^{a};q)_{k+\ell}}{(q^{d};q)_{k}(q;q)_{k}(p;p)_{\ell}}x^{k}y^{\ell}\int_{0}^{1}t^{b+\ell-1}\frac{(pt;p)_{\infty}}{(tp^{c-b};p)_{\infty}}d_{p}t\\
&=\frac{\Gamma_{p}(c)}{\Gamma_{p}(b)\Gamma_{p}(c-b)}\int_{0}^{1}\sum_{\ell,k=0}^{\infty}\frac{(q^{a};q)_{k+\ell}(p^{e};p)_{\ell}}{(p^{e};p)_{\ell}(q^{d};q)_{k}(q;q)_{k}(p;p)_{\ell}}x^{k}(yt)^{\ell}t^{b-1}\frac{(pt;p)_{\infty}}{(tp^{c-b};p)_{\infty}}d_{p}t\\
&=\frac{\Gamma_{p}(c)}{\Gamma_{p}(b)\Gamma_{p}(c-b)}\int_{0}^{1}t^{b-1}\Psi_{1}(q^{a},p^{e};p^{e},q^{d};q,p,x,yt)\frac{(pt;p)_{\infty}}{(tp^{c-b};p)_{\infty}}d_{p}t.
\end{split}
\end{eqnarray*}
Similarly, we obtain (\ref{2.39}).

Using the $p$-shifted factorials $(p^{b};p)_{\ell}$ and the $p$-Gamma functions are defined as follows (see \cite{kc, kk})
\begin{eqnarray*}
\begin{split}
(p^{b};p)_{\ell}=\frac{(1-p)^{\ell}\Gamma_{p}(b+\ell)}{\Gamma_{p}(b)}
\end{split}
\end{eqnarray*}
and
\begin{eqnarray*}
\begin{split}
\Gamma_{p}(b)=\int_{0}^{\frac{1}{1-p}}\;E_{p}(-pt)t^{b-1}d_{p}t,
\end{split}
\end{eqnarray*}
where $E_{p}(t)$ is the $q$-analogues of the exponential functions by
\begin{eqnarray*}
\begin{split}
E_{p}(t)=\sum_{r=0}^{\infty}q^{\frac{r(r-1)}{2}}\frac{t^{r}}{[r]_{q}!},
\end{split}
\end{eqnarray*}
we obtain (\ref{2.40})-(\ref{2.41}).
\end{proof}
\begin{thm} The summation formulas for $\Psi_{1}$ and $\Psi_{2}$ hold true
\begin{eqnarray}
\begin{split}
\Psi_{1}(q^{a},p^{b};p^{c},q^{d};q,p,x,y)=\frac{1}{(x;q)_{\infty}}\sum_{\ell=0}^{\infty}\frac{(q^{a};q)_{\ell}(p^{b};p)_{\ell}}{(p^{c};p)_{\ell}(p;p)_{\ell}}y^{\ell}\;_{1}\phi_{1}(q^{d-a-\ell};q^{d};q,xq^{a+\ell}),\label{2.42}
\end{split}
\end{eqnarray}
\begin{eqnarray}
\begin{split}
\Psi_{1}(q^{a},p^{b};p^{c},q^{d};q,p,x,y)=&\frac{1}{(q^{d};q)_{\infty}(x;q)_{\infty}}\sum_{\ell=0}^{\infty}\frac{(q^{a};q)_{\ell}(p^{b};p)_{\ell}}{(p^{c};p)_{\ell}(p;p)_{\ell}}(xq^{a+\ell};q)_{\infty}\\
&\times y^{\ell}\;_{1}\phi_{1}(x;xq^{a+\ell};q,q^{d}),\label{2.43}
\end{split}
\end{eqnarray}
\begin{eqnarray}
\begin{split}
\Psi_{2}(q^{a};p^{b},q^{c};q,p,x,y)=\frac{1}{(x;q)_{\infty}}\sum_{\ell=0}^{\infty}\frac{(q^{a};q)_{\ell}}{(p^{b};p)_{\ell}(p;p)_{\ell}}y^{\ell}\;_{1}\phi_{1}(q^{c-a-\ell};q^{c};q,xq^{a+\ell})\label{2.44}
\end{split}
\end{eqnarray}
and
\begin{eqnarray}
\begin{split}
\Psi_{2}(q^{a};p^{b},q^{c};q,p,x,y)=&\frac{1}{(q^{c};q)_{\infty}(x;q)_{\infty}}\sum_{\ell=0}^{\infty}\frac{(q^{a};q)_{\ell}}{(p^{b};p)_{\ell}(p;p)_{\ell}}(xq^{a+\ell};q)_{\infty}\\
&\times y^{\ell}\;_{1}\phi_{1}(x;xq^{a+\ell};q,q^{c}).\label{2.45}
\end{split}
\end{eqnarray}
\end{thm}
\begin{proof}
The $q$-extension of kummer's transformation formulas for the basic confluent hypergeometric function is obtained (see \cite{gr2, ks})
\begin{eqnarray}
\begin{split}
\;_{2}\phi_{1}(q^{a},0;q^{c};q,x)=\frac{1}{(x;q)_{\infty}}\;_{1}\phi_{1}(q^{c-a};q^{c};q,xq^{a})\label{2.46}
\end{split}
\end{eqnarray}
and
\begin{eqnarray}
\begin{split}
\;_{2}\phi_{1}(q^{a},0;q^{c};q,x)=\frac{(xq^{a};q)_{\infty}}{(q^{c};q)_{\infty}(x;q)_{\infty}}\;_{1}\phi_{1}(x;xq^{a};q,q^{c}).\label{2.47}
\end{split}
\end{eqnarray}
Using (\ref{2.31}) and (\ref{2.46}), we obtain (\ref{2.42}). Similarly, by using (\ref{2.32}), (\ref{2.46}) and (\ref{2.47}), we can prove (\ref{2.43})-(\ref{2.45}).
\end{proof}
\begin{rem}
The basic functions $\Psi_{1}$ and $\Psi_{2}$ are a $q$-analogue of the Humbert hypergeometric functions $\Psi_{1}$ and $\Psi_{2}$ defined by (\cite{emot1}, page 225, Equations (23)-(24)):
\begin{eqnarray}
\begin{split}
\lim_{q\longrightarrow 1}\Psi_{1}(q^{a},p^{b};p^{c},q^{d};q,p,x,\frac{1}{1-q}y)=&\Psi_{1}(a,b;c,d;x,y)\label{2.48}
\end{split}
\end{eqnarray}
and
\begin{eqnarray}
\begin{split}
\lim_{q\longrightarrow 1}\Psi_{2}(q^{a};p^{b},q^{c};q,p,x,\frac{1-p}{1-q}y)=&\Psi_{2}(a,b;c,d;x,y).\label{2.49}
\end{split}
\end{eqnarray}
\end{rem}
\begin{proof}
Using the limit 
\begin{equation*}
\begin{split}
\lim_{q\longrightarrow 1}\frac{(q^{a};q)_{n}}{(1-q)^{n}}&=(a)_{n},
\end{split}
\end{equation*}
we obtain (\ref{2.48})-(\ref{2.49}).
\end{proof}
\begin{thm}
The relationship holds true between the bibasic functions $\Psi_{1}$ and $\Psi_{2}$:
\begin{eqnarray}
\begin{split}
&\lim_{b\longrightarrow \infty}\Psi_{1}(q^{a},p^{b};p^{c},q^{d};q,p,x,y)&=\Psi_{2}(q^{a};p^{c},q^{d};q,p,x,y).\label{2.50}
\end{split}
\end{eqnarray}
\end{thm}
\begin{proof}
Using the limit
\begin{equation*}
\begin{split}
\lim_{b\longrightarrow \infty}(p^{b};p)_{n}&=(0;p)_{n}=1,
\end{split}
\end{equation*}
we obtain (\ref{2.50}).
\end{proof}
\subsection{Particular cases}
Here, we discuss and give theorems which are the following particular cases of theorem 2.10:
\begin{thm} The relationships between for $\Psi_{1}$, $\Psi_{2}$ and basic hypergeometric functions hold true
\begin{eqnarray}
\begin{split}
\Psi_{1}(q^{a},q^{b};q^{c},q^{d};q,0,y)=\;_{2}\phi_{1}(q^{a},q^{b};q^{c};q,y),\label{2.51}
\end{split}
\end{eqnarray}
\begin{eqnarray}
\begin{split}
\Psi_{1}(q^{a},q^{b};q^{c},q^{d};q,x,0)=\;_{2}\phi_{1}(q^{a},0;q^{d};q,x)\label{2.52}
\end{split}
\end{eqnarray}
and
\begin{eqnarray}
\begin{split}
\Psi_{2}(q^{a};q^{b},q^{c};q,x,0)=\;_{1}\phi_{0}(q^{a};-;q,x),\\
\Psi_{2}(q^{a};q^{b},q^{c};q,0,y)=\;_{2}\phi_{1}(q^{a},0;q^{b};q,y).\label{2.53}
\end{split}
\end{eqnarray}
\end{thm}
\begin{proof}
Set $p=q$ in (\ref{2.31}), we obtain
\begin{eqnarray*}
\begin{split}
\Psi_{1}(q^{a},q^{b};q^{c},q^{d};q,0,y)=\;_{2}\phi_{1}(q^{a},q^{b};q^{c};q,y).
\end{split}
\end{eqnarray*}
The identities (\ref{2.52})-(\ref{2.53}) can be proved in a like manner.
\end{proof}
\begin{thm} The summation formulas for $\Psi_{1}$ and $\Psi_{2}$ hold true:
\begin{eqnarray}
\begin{split}
\Psi_{1}(q^{a},q^{b};q^{c},q^{d};q,x,y)=&\frac{(q^{a};q)_{\infty}(q^{b};q)_{\infty}(q^{a}y;q)_{\infty}}{(q^{c};q)_{\infty}(q^{d};q)_{\infty}{(y;q)_{\infty}}}\sum_{r,s=0}^{\infty}\frac{(q^{d-a};q)_{r}(q^{c-b};q)_{s}(y;q)_{s}}{(q^{a}y;q)_{s}(q;q)_{s}(q;q)_{r}}\\
&\times q^{ar}q^{bs}\;_{2}\phi_{1}(0,0;q^{a+s}y;q,q^{r}x),\label{2.54}
\end{split}
\end{eqnarray}
\begin{eqnarray}
\begin{split}
\Psi_{1}(q^{a},q^{b};q^{c},q^{d};q,x,y)=&\frac{(q^{a};q)_{\infty}(q^{b};q)_{\infty}(xq^{a};q)_{\infty}}{(q^{c};q)_{\infty}(q^{d};q)_{\infty}(x;q)_{\infty}}\sum_{r,s=0}^{\infty}\frac{(q^{d-a};q)_{r}(q^{c-b};q)_{s}}{(q;q)_{s}(q;q)_{r}}\\
&\times q^{ar}q^{bs}\;_{2}\phi_{1}(0,0;xq^{a};q,q^{r}q^{s}y),\\
\Psi_{1}(q^{a},q^{b};q^{c},q^{d};q,x,y)=&\frac{(q^{a};q)_{\infty}(q^{b};q)_{\infty}}{(q^{c};q)_{\infty}(q^{d};q)_{\infty}}\sum_{r,s,\ell=0}^{\infty}\frac{(q^{d-a};q)_{r}(q^{c-b};q)_{s}}{(q;q)_{s}(q;q)_{r}(q;q)_{\ell}}\\
&\times q^{ar}q^{bs}\bigg{(}q^{s}y\bigg{)}^{\ell}\:_{1}\Phi_{0}(q^{a+\ell};-;q,q^{\ell}x)\label{2.55}
\end{split}
\end{eqnarray}
and
\begin{eqnarray}
\begin{split}
&\Psi_{2}(q^{a};q^{b},q^{c};q,x,y)=\frac{(q^{a};q)_{\infty}}{(q^{b};q)_{\infty}}\sum_{\ell=0}^{\infty}\frac{(q^{b+\ell};q)_{\infty}}{(q^{a+\ell};q)_{\infty}(q;q)_{\ell}}y^{\ell}\;_{2}\phi_{1}(q^{a+\ell},0;q^{c};q,x),\\
&\Psi_{2}(q^{a};q^{b},q^{c};q,x,y)=\frac{(q^{a};q)_{\infty}}{(q^{c};q)_{\infty}}\sum_{k=0}^{\infty}\frac{(q^{c+k};q)_{\infty}}{(q^{a+k};q)_{\infty}(q;q)_{k}}x^{k}\;_{2}\phi_{1}(q^{a+k},0;q^{b};q,y).\label{2.56}
\end{split}
\end{eqnarray}
\end{thm}
\begin{proof}
Take $p=q$ with a similar argument as in theorem 2.12 and use to transform the series on the right, simplification gives the required results (\ref{2.54})-(\ref{2.56}).
\end{proof}
\begin{cor} The summation formula for $\Psi_{1}$ hold true
\begin{eqnarray}
\begin{split}
\Psi_{1}(q^{a},q^{b};q^{c},q^{d};q,x,y)=&\frac{(q^{a};q)_{\infty}(q^{b};q)_{\infty}(q^{a}y;q)_{\infty}}{(q^{d};q)_{\infty}(q^{c};q)_{\infty}{(y;q)_{\infty}}}\sum_{r,s=0}^{\infty}\frac{(q^{d-a};q)_{r}(q^{c-b};q)_{s}(y;q)_{s}}{(q^{a}y;q)_{s}(q;q)_{s}(q;q)_{r}}\\
&\times \frac{q^{ar+bs}}{(q^{a+s}y;q)_{\infty}(xq^{r};q)_{\infty}}\;_{1}\phi_{1}(xq^{r};0;q,q^{a+s}y),\\
\Psi_{1}(q^{a},q^{b};q^{c},q^{d};q,x,y)=&\frac{(q^{a};q)_{\infty}(q^{b};q)_{\infty}(q^{a}y;q)_{\infty}}{(q^{d};q)_{\infty}(q^{c};q)_{\infty}{(y;q)_{\infty}}}\sum_{r,s=0}^{\infty}\frac{(q^{d-a};q)_{r}(q^{c-b};q)_{s}(y;q)_{s}}{(q^{a}y;q)_{s}(q;q)_{s}(q;q)_{r}}\\
&\times \frac{q^{ar+bs}}{(xq^{r};q)_{\infty}}\;_{0}\phi_{1}(-;q^{a+s}y;q,xq^{r}q^{a+s}y),\label{2.57}
\end{split}
\end{eqnarray}
and
\begin{eqnarray}
\begin{split}
\Psi_{1}(q^{a},q^{b};q^{c},q^{d};q,x,y)=&\frac{(q^{a};q)_{\infty}(q^{b};q)_{\infty}(xq^{a};q)_{\infty}}{(q^{d};q)_{\infty}(q^{c};q)_{\infty}(x;q)_{\infty}}\sum_{r,s=0}^{\infty}\frac{(q^{d-a};q)_{r}(q^{c-b};q)_{s}}{(q;q)_{s}(q;q)_{r}}\\
&\times \frac{q^{ar+bs}}{(xq^{a};q)_{\infty}(q^{r+s}y;q)_{\infty}}\;_{1}\phi_{1}(q^{r+s}y;0;q,xq^{a}),\\
\Psi_{1}(q^{a},q^{b};q^{c},q^{d};q,x,y)=&\frac{(q^{a};q)_{\infty}(q^{b};q)_{\infty}(xq^{a};q)_{\infty}}{(q^{d};q)_{\infty}(q^{c};q)_{\infty}(x;q)_{\infty}}\sum_{r,s=0}^{\infty}\frac{(q^{d-a};q)_{r}(q^{c-b};q)_{s}}{(q;q)_{s}(q;q)_{r}}\\
&\times \frac{q^{ar+bs}}{(q^{r+s}y;q)_{\infty}}\;_{0}\phi_{1}(-;xq^{a};q,xyq^{a+r+s}).\label{2.58}
\end{split}
\end{eqnarray}

\end{cor}
\begin{proof}
Using the Heine's transformations formulas for basic hypergeometric functions:(see \cite{gr2})
\begin{eqnarray}
\begin{split}
\;_{2}\phi_{1}(0,0;q^{c};q,x)&=\frac{1}{(q^{c};q)_{\infty}(x;q)_{\infty}}\;_{1}\phi_{1}(x;0;q,q^{c}),\\
&=\frac{1}{(x;q)_{\infty}}\;_{0}\phi_{1}(-;q^{c};q,xq^{c}).\label{2.59}
\end{split}
\end{eqnarray}
and replacing $q^{c}$ and $x$ by $q^{a+s}y$ and $xq^{r}$ in (\ref{2.59}), we get (\ref{2.57}), replacing $q^{c}$ and $x$ by $xq^{a}$ and $q^{r+s}y$ in (\ref{2.59}), we obtain (\ref{2.58}).
\end{proof}
%
\section{Concluding remarks}
It may be remarked, finally, by using the same technique in above By new definitions of bibasic Humbert hypergeometric functions $\Psi_{1}$ and $\Psi_{2}$ for the numerator and denominator parameters with a few more new and elegant ones on quantum analogue. Furthermore, we obtained interesting results. We believe that a detailed study of the bibasic Humbert hypergeometric functions on properties for this functions should be very interesting to permit elegant and neat extensions and hence have not been included herein.


\end{document}